\documentclass[11pt]{article}

\usepackage[T1]{fontenc}
\usepackage[latin1]{inputenc}
\usepackage{mathtools,amsthm,amssymb}

\usepackage{mathrsfs}

\usepackage{authblk}
\usepackage{xcolor}

\usepackage{hyperref}
\usepackage[margin=1in]{geometry}

\newtheorem{hypothesis}{Hypothesis}
\newtheorem*{hypothesisRemark}{Hypothesis-Remark}

\newtheorem{theo}{Theorem}[section]
\newtheorem{cor}[theo]{Corollary}
\newtheorem{lem}[theo]{Lemma}

\def \PP{\mathbb{P}}
\def \RR{\mathbb{R}}
\def \EE{\mathbb{E}}
\newcommand{\Ea}{ {\cal E }}
\newcommand{\Sa}{ {\cal S }}
\newcommand{\Fa}{ {\cal F }}
\newcommand{\Ma}{ {\cal M }}
\newcommand{\Pa}{ {\cal P }}
\newcommand{\Aa}{ {\cal A }}
\newcommand{\Wa}{ {\cal W }}
\newcommand{\Ra}{ {\cal R }}
\newcommand{\Ua}{ {\cal U }}
\newcommand{\Va}{ {\cal V }}
\newcommand{\Ga}{ {\cal G }}
\def\tr{\mbox{\rm tr}}
\newcommand{\Ca}{ {\cal C }}

\begin{document}

\title{Stability Properties of Systems of Linear Stochastic Differential Equations with Random Coefficients}

\author[$1$]{Adrian N. Bishop}
\author[$2$]{Pierre Del Moral}
\affil[$1$]{{\small University of Technology Sydney (UTS); and CSIRO, Australia}}
\affil[$2$]{{\small INRIA, Bordeaux Research Center, France; and UNSW, Sydney, Australia}}
\date{}

\maketitle

\begin{abstract}
This work is concerned with the stability properties of linear stochastic differential equations with random (drift and diffusion) coefficient matrices, and the stability of a corresponding random transition matrix (or exponential semigroup). We consider a class of random matrix drift coefficients that involves random perturbations of an exponentially stable flow of deterministic (time-varying) drift matrices. In contrast with more conventional studies, our analysis is not based on the existence of Lyapunov functions, and it does not rely on any ergodic properties. These approaches are often difficult to apply in practice when the drift/diffusion coefficients are random. We present rather weak and easily checked perturbation-type conditions for the asymptotic stability of time-varying and random linear stochastic differential equations. We provide new log-Lyapunov estimates and exponential contraction inequalities on any time horizon as soon as the fluctuation parameter is sufficiently small. These seem to be the first results of this type for this class of linear stochastic differential equations with random coefficient matrices. 
\end{abstract}

\section{Introduction}\label{intro-sec}
\subsection{Aims and Scope}
This work is concerned with the stability of a system of linear stochastic differential equations (i.e. Ito-type diffusion equations driven by a Wiener process), with the added complexity of time-varying random coefficients. That is, we consider the stability of general systems of the form
\begin{equation}\label{intro-X}
	dX_t^{\epsilon}=A_t^{\epsilon}\,X_t^{\epsilon}~dt+(B_t^{\epsilon})^{1/2}\,dW_t
\end{equation}
where $W_t$ an $r$-dimensional Wiener process and $(A_t^{\epsilon},B^{\epsilon}_t)$ are suitably defined (non-anticipating) random processes indexed by some parameter $\epsilon\in[0,1]$ and independent of $(X_0^{\epsilon},W_t)$. Detailed technical models are given later. We term these processes: random\footnote{The ``random'' predicate is used to denote the randomness of the coefficients. In the case of deterministic coefficients, we term this equation: a stochastic differential equation, or an Ito-type (linear) diffusion process (driven by a Wiener process). Note that Ito-type integrals do not preclude (suitably adapted) random integrands; however, our terminology is used to distinguish those equations with random coefficients and those without.} linear stochastic differential equations, or random linear diffusions, or random Ornstein-Uhlenbeck processes\footnote{``Random'' is again used to describe the randomness of the coefficients and we refer to this process as an Ornstein-Uhlenbeck process since we are concerned ultimately with stable processes in an Ornstein-Uhlenbeck-type form.}. The regularity and stability properties of this process are far from obvious. 

The stability properties of Ito-type diffusion equations with deterministic coefficients (i.e. Wiener processes with non-zero deterministic drift) has been well studied, and we point to the texts \cite{kushner1967book,Kozin69,khasminskii2011stochastic,mao2010stochastic} for a detailed survey of results. Rather independently, systems of ordinary differential equations with random coefficients (i.e. systems of the form above with $B_t^{\epsilon}=0$) have also been considered, and we point to \cite{Bertram1959,kats60,Infante68,Kozin69,Kozin73,Tsokos74,Blankenship77,Geman79,solo69} for a collection of results and techniques. Stability results for these latter-type of equations are often difficult to apply in practice. Stability of the special case concerning ordinary differential equations with piecewise constant (random) coefficients was studied in, e.g., \cite{Bharucha61,Morozan67} and is applicable to stochastic jump systems \cite{fang2002}. These reference lists are not exhaustive (see also the references therein). We return to some of this literature later, with a more specific technical relationship to our current work. 

Combined diffusion-type equations with random coefficients are more general, and stability results are difficult to apply in practice (owing a lot again to the randomness in the coefficients). Nevertheless, these types of model appear (in some fashion) in economics and finance \cite{Turnovsky73,lim02}, biology \cite{Tsokos74,Tiwari76}, mechanics and physics \cite{Soong73}, etc.

These models are relevant also in mathematical control and systems theory, and in estimation, filtering and data assimilation. Early work by Bismut \cite{bismut1976} considered the linear-quadratic optimal control of very general (linear) systems of this form, and in which the control may act on both the drift and the diffusion. Optimal control in this framework is complicated by the need to address the well-posedness of certain backward stochastic differential equations with (the added complexity of) random coefficients. Much work has been considered in this direction and it is mostly beyond the scope of this article; we point to \cite{yong1999stochastic,Tang2005,Tang2015,Pham2016} and the references therein.  

Nonlinear analogues of (\ref{intro-X})  (i.e. nonlinear stochastic differential equations driven by a Wiener process with random drift/diffusion functions or random inputs to the drift/diffusion) are also of interest, e.g. \cite{Soong73,Tsokos74,peng92,yong1999stochastic}, albeit they are beyond the scope of this work.

Infinite dimensional analogues of (\ref{intro-X}) (e.g. linear diffusions taking values in Banach spaces) are related to stochastic partial differential equations. In this context, the diffusion is defined in terms of a cylindrical Brownian motion on some Hilbert space and the drift matrix is replaced by some generator of a strongly continuous semigroup on some Banach space~\cite{daprato-zabczyk,daprato-zabczyk-2}. The stability properties of this class of stochastic partial differential equations is technically much more involved and they are also beyond the scope of this work.

The preceding literature review is primarily focused on general stability results concerning random linear stochastic differential equations, or random linear diffusions. 

\subsection{An Illustrative Example: Ensemble Kalman-Bucy Filters}

Our main results are generally applicable to systems of the form (\ref{intro-X}). We illustrate the model (\ref{intro-X}), and our regularity conditions (given later), via a running example in filtering and data assimilation; i.e. by studying the stability properties of ensemble Kalman-Bucy-type filters \cite{evensen03}. See \cite{apa-2017,2017arXiv171110065B,2018arXiv180800235} for dedicated technical discussions on related ensemble filters, and associated literature review (which is beyond the focus of this article). Related applications of our type of stability result are in \cite{Guo94}. 

The ensemble Kalman-Bucy filter \cite{evensen03} (in the linear-Gaussian problem setting) is a mean-field particle approximation of the classical Kalman-Bucy filter; see \cite{DelMoral/Tugaut:2016,Bishop/DelMoral:2016}. This approximation involves a collection of $N\geq2$ interacting particles driven by Kalman-Bucy-type update diffusions, but with an interaction function (i.e.~the Kalman gain) that depends on the (random) sample covariance. 
 
 More precisely, up to a change of probability space, the $N$-sample covariances matrices associated with these filters
 satisfy a nonlinear matrix-valued (Riccati-type) diffusion equation of the form,
  \begin{equation}\label{ric-intro}
d\Pa^{\epsilon}_t=\left[\Aa\,\Pa^{\epsilon}_t+\Pa^{\epsilon}_t\,\Aa^{\prime}+\Ra-\Pa^{\epsilon}_t\,\Sa\,\Pa^{\epsilon}_t\right]dt+\epsilon\left[
\left(\Pa^{\epsilon}_t\right)^{1/2}\,d\Wa_t~\Sigma^{1/2}\left(\Pa^{\epsilon}_t\right)+\Sigma^{1/2}\left(\Pa^{\epsilon}_t\right)d\Wa_t^{\prime}\left(\Pa^{\epsilon}_t\right)^{1/2}\right] \end{equation}
with a fluctuation parameter $\epsilon=\epsilon_N\longrightarrow_{N\rightarrow\infty} 0$ and $(r\times r)$-matrices $\Wa_t$ with independent Brownian entries. Here, $(\cdot)^{\prime}$ denotes the transposition operator, $(\Aa,\Ra,\Sa)$ are $(r\times r)$-matrices, $\Ra,\Sa>0$ are positive definite matrices, and $\Sigma(\cdot)$ is a positive map from the set of symmetric $(r\times r)$-matrices into the set of positive semi-definite $(r\times r)$-matrices such that,
\begin{equation}\label{sigma-intro}
	0 \,\leq\, \Sigma(\Pa) \,\leq\, \Ua+\Pa\,\Va\,\Pa\quad \mbox{\rm for some positive semi-definite matrices}\quad \Ua,\Va \,\geq\, 0
\end{equation}
The matrices $\Aa$ and $\Ra$ capture the drift and the covariance diffusion matrices of the underlying signal. The matrix $\Sa:=\Ca^{\prime}\Ra_o^{-1}\,\Ca$ is defined in terms of the sensor observation matrix $\Ca$, and the covariance diffusion matrix $\Ra_o$ of the observation noise. The dominating matrices $(\Ua,\Va)$ depend on the different variants of ensemble Kalman-Bucy filters one may consider~\cite{2018arXiv180800235}.

    Now, the difference between the ensemble Kalman-Bucy filter sample mean and the true signal state (i.e. the estimation error) is described by a stochastic process $X^{\epsilon}_t$ of the form (\ref{intro-X}), for some Brownian motion $W_t$ independent of $\Wa_t$, and with the drift and diffusion matrices defined by
  \begin{equation}\label{KB-intro}
A_t^{\epsilon}=\Aa-\Pa^{\epsilon}_t\,\Sa\quad\mbox{\rm and}\quad B_t^{\epsilon}=\Sigma_{\epsilon}(\Pa^{\epsilon}_t)
\end{equation}
Here, $\Sigma_{\epsilon}$ denotes another (related) positive map satisfying (\ref{sigma-intro}) for some matrices $\Ua_{\epsilon},\Va_{\epsilon}$ that are uniformly bounded w.r.t. the parameter $\epsilon$.

To relate this derivation to classical Kalman-Bucy filtering, note if $\epsilon=0$, then the diffusion matrix $\Pa_t:=\Pa^{0}_t$ reduces to the matrix Riccati differential equation,
\begin{equation}\label{intro-true-Ricc}
 \partial_t\,\Pa_t=\mbox{\rm Ricc}(\Pa_t):=\Aa\,\Pa_t+\Pa_t\,\Aa^{\prime}+\Ra-\Pa_t\,\Sa\,\Pa_t
\end{equation}  
Moreover, the stochastic process $X_t:=X^{0}_t$ defined in (\ref{intro-X}) with the matrices (\ref{KB-intro}) and $\epsilon=0$ is just the difference between the classical Kalman-Bucy state estimate (i.e. the conditional expectation of the state given the observations) and the true signal state. In this context, note that $\Pa_t$ is exactly the covariance matrix of the error process $X_t$. When the pair of matrices $(\Aa,\Ra^{1/2})$ is stabilisable, and $(\Aa,\Sa^{1/2})$ is detectable, then the error $X_t$ is a stable process in the sense that the drift matrix
\begin{equation}\label{ref-intro-H0-0}
	A_t:=\Aa-\Pa_t\,\Sa
\end{equation}  
delivers a uniformly exponentially stable linear (time-varying) system. Moreover, there exists a unique matrix $\Pa_{\infty}\geq 0$ such that
\begin{equation}\label{ref-intro-H0-1}
\mbox{\rm Ricc}(\Pa_{\infty})=0\quad \mbox{\rm and the spectral abscissa of the matrix ~~$A_{\infty}:=\Aa-\Pa_{\infty}\,\Sa$~~ is negative}
\end{equation} 
If $(\Aa,\Ra^{1/2})$ is controllable, then $\Pa_{\infty}>0$. When the pair $(\Aa,\Ra^{1/2})$ is stabilisable, and  $(\Aa,\Sa^{1/2})$ is detectable, 
we have
\begin{equation}\label{ref-intro-H0-2}
\mbox{\rm $\Pa_t~\longrightarrow_{t\rightarrow\infty}~\Pa_{\infty}$~ and thus ~$A_t~\longrightarrow_{t\rightarrow\infty}~A_{\infty}$~
exponentially fast for any initial matrices}
\end{equation} 
For proof of these classical results on the true Kalman-Bucy filter see, e.g., \cite{kucera72,Lancaster1995}, and the convergence results in \cite{Kwakernaak72,callier81,vanhandel2009,Bishop/DelMoral:2016,bd-CARE}. 

Going forward we may assume an even stronger notion of controllability and observability; and suppose that the logarithmic norm (defined later) of $A_{\infty}$ is negative. A simple example of this type is to suppose that there exists a change of basis such that $\Aa$ is symmetric and $\Sa\propto \mathrm{Id}$.

 Returning to the mean-field approximation and the ensemble Kalman-Bucy filter. The Riccati diffusion equation (\ref{ric-intro}) encapsulates several classes of ensemble Kalman filter, including the continuous time versions of the vanilla ensemble Kalman filter~\cite{evensen03} as well as the so-called deterministic ensemble Kalman filter introduced in~\cite{sakov2008a}. It also applies to various classes of regularized ensemble filter discussed in~\cite{Bishop/DelMoral/Pathiraja:2017}; e.g. so-called covariance inflation-type regularisation  \cite{anderson99,hamill01,evensen03}.  
 
  Under further natural observability conditions (e.g. $\Sa\propto \mathrm{Id}$), the sample covariance in (\ref{ric-intro}) can be made as close as desired to the solution of the stable Riccati equation (\ref{intro-true-Ricc}); see~\cite{apa-2017,2018arXiv180800235,DelMoral/Tugaut:2016}. One advantage of the deterministic-type ensemble Kalman filters \cite{sakov2008a}, is that $\Va=0$; so the sample covariance $\Pa^{\epsilon}_t$ exhibits less fluctuation around its limiting value $\Pa_t$ than in the vanilla case \cite{evensen03}. For details on these models and their fluctuation properties see~\cite{apa-2017,2017arXiv171110065B,2018arXiv180800235} and the references therein.

As an illustrative example of our general stability results, we will relate our analysis of random linear diffusions to those associated with the ensemble Kalman-Bucy filter, and the drift and diffusion matrices (\ref{KB-intro}). Related analysis in one-dimension is given in \cite{2017arXiv171110065B} where stronger results are available.

\subsection{General Models and Notation}

We denote by $\Ma_{r}=\RR^{r\times r}$ the set of $(r\times r)$-square matrices with real entries and $r\geq 1$. We let $\Sa_r\subset \Ma_{r}$ denote the subset of symmetric matrices, $\Sa_r^0\subset\Sa_r$ the subset of positive semi-definite matrices, and $\Sa_r^+\subset \Sa_r^0$ the subset of positive definite matrices.  Given $B\in \Sa_r^0-\Sa_r^+$ we denote by $B^{1/2}$ a (non-unique) but symmetric square root of $B$ (given by a Cholesky decomposition). When $B\in\Sa_r^+$ we always choose the principal (unique) symmetric square root. We write $A^{\prime}$ to denote the transposition of a matrix $A$, and $A_{sym}=(A+A^{\prime})/2$ to denote the symmetric part of $A\in\Ma_{r}$.

We equip the set $\Ma_{r}$ with the spectral norm $\Vert A \Vert:=\Vert A \Vert_2=\sqrt{\lambda_{1}(AA^{\prime})}$ where $\lambda_{1}(\cdot)$ denotes the maximal eigenvalue. Let $\tr(A)=\sum_{1\leq i\leq r}A(i,i)$ denote the trace operator. We also denote by $\mu(A)=\lambda_{1}(A_{sym})$ its logarithmic norm. Note that $$\mu(A)~\geq~\max_i{\left\{\mbox{\rm Re}[\lambda_i(A)]\right\}}$$ 

Throughout the remainder, $A~:~t\in \RR_+:=[0,\infty[~\mapsto A_t\in \Ma_r$ denotes some deterministic flow of matrices satisfying  the following condition:
\begin{hypothesis}[$H_0$]
$$
\qquad \Vert A_t-A_{\infty}\Vert ~\leq~ a\,e^{-bt}\quad\mbox{\rm for some}~A_{\infty}\in \Ma_r\quad\mbox{s.t.}\quad
\mu(A_{\infty})<0
$$
for any time horizon $t\geq 0$ and some parameters $a<\infty$ and $b>0$. We set $c_0:=\vert \mu(A_{\infty})\vert$ throughout. 
\end{hypothesis}

For example, the drift matrices $A_t$ of the true Kalman-Bucy error process defined in (\ref{ref-intro-H0-0}) satisfy the above condition under suitably strong controllability and observability conditions, and with the limiting matrix, $\epsilon=0$, $t\rightarrow\infty$, denoted by $A_{\infty}$ and defined in (\ref{ref-intro-H0-1}).

Let $\Ea_{s,t}(A)$ be the exponential semigroup (or the state transition matrix) associated with a smooth flow of matrices $A:t\in\RR_+\mapsto A_t\in \Ma_r$ defined for any $s\leq t$ by the forward and backward differential equations,
\begin{equation*}
 \partial_t \,\Ea_{s,t}(A)=A_t\,\Ea_{s,t}(A)\quad\mbox{\rm and}\quad
\partial_s\, \Ea_{s,t}(A)=-\Ea_{s,t}(A)\,A_s
\end{equation*}
with $\Ea_{s,s}(A)=\mathrm{Id}$, the identity matrix (of appropriate dimension). Equivalently in terms of the matrices
$\Ea_t(A):=\Ea_{0,t}(A)$ we have
$
\Ea_{s,t}(A)=\Ea_t(A)\Ea_s(A)^{-1}
$. 

For any $s,t>0$ we recall the logarithmic norm estimate
\begin{equation}
\log{\Vert \Ea_{s,s+t}(A)\Vert}\leq \int_s^{s+t}\mu(A_u)~du~~\Longrightarrow~~
 \frac{1}{t}\log{\Vert \Ea_{s,s+t}(A)\Vert}\leq \mu(A_{\infty})+\frac{e^{-bs}}{t}~\frac{a}{b} \label{stabilityestimate1}
\end{equation}
The l.h.s. log-norm estimate comes from the fact that
$$
\partial_tx_t:=A_tx_t\quad\Longrightarrow \quad\frac{1}{2}\,\partial_t\,\Vert x_t\Vert^2~=~x^{\prime}_t(A_t)_{sym}x_t~\leq ~\mu(A_t)~\Vert x_t\Vert^2
$$
In the above display, $\Vert x\Vert$ denotes the Euclidian norm of a vector $x\in\RR^r$. The extension of the log-norm estimate (\ref{stabilityestimate1}) to any norm on $\RR^r$ and any matrix norm can be found in~\cite{coppel1975}.
\begin{hypothesisRemark}[$H_0'$]
With the l.h.s. of this implication (\ref{stabilityestimate1}) in mind, a straightforward extension would be to consider a relaxing of $(H_0)$ to the case in which the flow associated with $A_t$ has no time-limiting fixed point, but is itself just a time-varying stabilising matrix. 
\end{hypothesisRemark}

\emph{We come to the random processes of interest.} Let $A^{\epsilon}:t\in\RR_+\mapsto A_t^{\epsilon}\in \Ma_r$ and $B^{\epsilon}:t\in\RR_+\mapsto B_t^{\epsilon}\in \Sa^0_r$ be a collection of c\`adl\`ag random processes defined on a common filtered probability space
$(\Omega,(\Fa_t^\epsilon)_{t\geq 0},\PP)$ and indexed by some parameter $\epsilon\in [0,1]$. Let $
\Fa_t^\epsilon:=\sigma((A^{\epsilon}_s,B^{\epsilon}_s),~s\leq t)$. We require,
 \begin{equation}\label{hyp-unif-moments}
\sup_{0\leq \epsilon\leq \epsilon_{0,n}}{~\sup_{t\geq 0}{\EE\left(\Vert A^{\epsilon}_t\Vert^n\right)}}<\infty\quad\mbox{\rm and}\quad
\rho_n:=\sup_{0\leq \epsilon\leq \epsilon_{1,n}}{~\sup_{t\geq 0}{\EE\left[\tr\left(B^{\epsilon}_t\right)^n\right]^{1/n}}}<\infty
 \end{equation}
for some fluctuation parameters $\epsilon_{0,n}$ and $\epsilon_{1,n}\in [0,1]$, and for any parameter $n\geq 1$. 

For example, in the context of ensemble Kalman-Bucy filters, the random matrices defined in (\ref{KB-intro}) satisfy the uniform estimates (\ref{hyp-unif-moments}) (cf. Theorem 2.2 in~\cite{2018arXiv180800235}, see also \cite{apa-2017}).

We let $X_t^{\epsilon}$ be the collection of  random Ornstein-Uhlenbeck process defined by
\begin{equation}\label{stochastic-OU}
dX_t^{\epsilon}=A_t^{\epsilon}~X_t^{\epsilon}~dt+(B_t^{\epsilon})^{1/2}~dW_t
\end{equation}
where $W_t$ an $r$-dimensional Wiener process and we assume that $(X_0^{\epsilon},W_t)$ are independent of the stochastic processes $(A_t^{\epsilon},B^{\epsilon}_t)$. We also denote by $X^{\epsilon,x}_t$ the process starting at $X^{\epsilon,x}_0=x\in\RR^r$, 

\emph{The objective of this work is to study the stability properties of the semigroup $\Ea_{s,t}(A^{\epsilon})$ associated with the stochastic process $A^{\epsilon}$, and the stability of the random Ornstein-Uhlenbeck process (\ref{stochastic-OU}).}

The analysis of the long time behaviour of the stochastic model (\ref{stochastic-OU}) differs strongly from the analysis of conventional time-invariant and deterministic linear dynamical systems. As with general time-varying deterministic linear dynamical systems, the asymptotic behaviour of (\ref{stochastic-OU}) cannot be characterised by the statistical properties of the spectral abscissa of the random matrices $A^{\epsilon}_t$. Indeed, unstable semigroups $\Ea_{s,t}(A)$ associated with time-varying (deterministic) matrices $A_t$ with negative eigenvalues are exemplified in~\cite{coppel1978stability,wu1974note}. Conversely, stable semigroups $\Ea_{s,t}(A)$ with $A_t$ having positive eigenvalues are given by Wu in~\cite{wu1974note}. The same general conclusion holds in the statistical case without (quite strong) additional restrictions on the class of model considered. We seek quite weak and more readily verifiable and practical conditions in this work.

Observe that the solution of (\ref{stochastic-OU}) is provided by the formula
\begin{equation}\label{solution-stochastic-OU}
X^{\epsilon}_t=\Ea_{t}(A^{\epsilon})X^{\epsilon}_0+X^{\epsilon,0}_t\quad \mbox{\rm with}\quad X^{\epsilon,0}_t=
\int_0^t~\Ea_{s,t}(A^{\epsilon})~(B^{\epsilon}_s)^{1/2}~dW_s
\end{equation}
The process $X^{\epsilon,0}_t$ defined above is sometimes called a stochastic convolution. Note that given $\Fa_t^\epsilon$, the r.h.s. integral in  (\ref{solution-stochastic-OU}) is an Ito integral since the Brownian motion $(W_s)_{s\leq t}$ is independent of $(A_s^{\epsilon},B^{\epsilon}_s)_{s\leq t}$; i.e. it is non-anticipative due to the independence of the relevant randomness. Nevertheless, the process $X^{\epsilon,0}_t$ in (\ref{solution-stochastic-OU}) is not a martingale (w.r.t. the filtration $\Ga_t$ generated by the Brownian motion, nor w.r.t. the enlarged filtration  $\Ga_t\vee\Fa_t^\epsilon$) and the analysis of its regularity properties is far from obvious.

Also note that
$$
\EE\left(X_t^{\epsilon}~|~\Fa_t^{\epsilon}\right)=\Ea_{t}(A^{\epsilon})\,\EE(X^{\epsilon}_0)
$$
This implies that the conditional covariance process is given by,
\begin{eqnarray*}
C_t^{\epsilon}&:=&\EE\left([X_t^{\epsilon}-\EE\left(X_t^{\epsilon}~|~\Fa_t^{\epsilon}\right)][X_t^{\epsilon}-\EE\left(X_t^{\epsilon}~|~\Fa_t^{\epsilon}\right)]^{\prime}~|~\Fa^{\epsilon}_t\right)\\
&=&  \Ea_{t}(A^{\epsilon})~C_0^{\epsilon}~\Ea_{t}(A^{\epsilon})^{\prime}+
\int_0^t\int_0^t~
\Ea_{s,t}(A^{\epsilon})~(B^{\epsilon}_s)^{1/2}~\underbrace{\EE(dW_s dW_r)}_{\mathbf{1}_{\{r=s\}}\,\mathrm{Id}\,ds}~(B^{\epsilon}_r)^{1/2}~\Ea_{r,t}(A^{\epsilon})^{\prime} \\
&=&  \Ea_{t}(A^{\epsilon})~C_0^{\epsilon}~\Ea_{t}(A^{\epsilon})^{\prime}+\int_0^t~\Ea_{s,t}(A^{\epsilon})~B^{\epsilon}_s~\Ea_{s,t}(A^{\epsilon})^{\prime}~ds
\end{eqnarray*}

In the special case when $B_t^{\epsilon}=0$, the Ornstein-Uhlenbeck process (\ref{stochastic-OU}) resumes to the linear random dynamical equation,
\begin{equation}\label{stochastic-OU-linear}
\partial_t X_t^{\epsilon}=A_t^{\epsilon}~X_t^{\epsilon}~\Longleftrightarrow~X_t^{\epsilon}=\Ea_{t}(A^{\epsilon})~X_0^{\epsilon}
\end{equation}
whose stability properties have been considered in, e.g., \cite{Bertram1959,kats60,Bharucha61,Morozan67,Infante68,Kozin69,Kozin73,Tsokos74,Blankenship77,Geman79}. As previously noted, stability results for systems of the form (\ref{stochastic-OU-linear}) are often difficult to apply in practice. It is also worth noting that the discrete-time version of (\ref{stochastic-OU-linear}) is given by the recursive equation
$$
X_{n}^{\epsilon}=A_n^{\epsilon}X_{n-1}^{\epsilon}=A_n^{\epsilon}A_{n-1}^{\epsilon}\ldots A_1^{\epsilon}X_{0}^{\epsilon}
$$
If $A_n^{\epsilon}$ is matrix-valued Markov chain, then for ergodic chains, the stability properties of this equation are related to Oseledec's multiplicative ergodic theorem~\cite{ludwig1998random,oseledec1968multiplicative}. More generally, the stability properties in discrete-time may be related to various theorems concerning the infinite product of stochastic matrices; e.g. see \cite{Wolfowitz63,moreau:05,saber:07}. See also \cite{Guo94} for related results in discrete-time.

\section{Statement of the Main Results}
\subsection{Some Regularity Conditions}

Coming to the stability properties of (\ref{stochastic-OU}), if we assume that for any $s\leq t$ and some $n\geq 4$ we have  the exponential estimate
\begin{equation}\label{ref-intro-perfect-est}
\EE\left(\Vert\Ea_{s,t}(A^{\epsilon})\Vert^n\right)^{1/n}\leq \alpha\,e^{-\beta (t-s)}
\end{equation} for some parameters $\alpha<\infty$ and $\beta>0$; then, for any $\epsilon\leq \epsilon_{1,2n}$ we find 
\begin{equation}\label{ref-intro-perfect-est-conseq}
\EE\left(\Vert X^{\epsilon,x}_t-X^{\epsilon,y}_t\Vert^n\right)^{1/n}\leq \alpha\, e^{-\beta (t-s)}~\Vert x-y\Vert\qquad \mbox{\rm and}\qquad
\sup_{t\geq 0}\EE\left(\Vert X^{\epsilon,0}_t\Vert^{n/2}\right)^{2/n}<\infty
\end{equation}
This elementary result indicates that the stability properties of the random Ornstein-Uhlenbeck process $X_t^{\epsilon}$ are directly related to the contraction properties of the stochastic semigroup $\Ea_{s,t}(A^{\epsilon})$ with $s\leq t$. This is of course not surprising. 

Nevertheless, in practice, exponential estimates of the form (\ref{ref-intro-perfect-est}) are difficult to obtain mainly because the semigroup $\Ea_{s,t}(A^{\epsilon})$ cannot be represented as elementary matrix exponentials, but rather only in terms of Peano-Baker series ~\cite{Peano} or sophisticated Magnus exponential series~\cite{blanes,magnus}; see also the studies~\cite{Brockett,Frazer,Ince}.  The analysis of the exponential semigroups associated with the stochastic matrices $A^{\epsilon}_t$ discussed in (\ref{KB-intro}) provides an intuitive feel of the complexity of these models.

We note that the l.h.s. uniform moment condition in (\ref{hyp-unif-moments})
ensures that the random drift process $ A^{\epsilon}_t$ is uniformly tight, in the sense that for any $\nu\in[0,1]$, $\exists k$ such that
$
\sup_{t\geq 0}{\PP\left(\Vert A^{\epsilon}_t\Vert\geq k \right)}\leq \nu 
$.
By Prohorov's theorem this implies that the distributions of the random states $(A^{\epsilon}_t)_{t\geq 0}$ is relatively compact so there exists at least one limiting invariant distribution $\pi_{\epsilon}$ on $\Ma_r$. In addition there exists a sequence of random times $t_n$ such that 
$
\mbox{\rm Law}(A^{\epsilon}_{t_n})~\longrightarrow_{n\rightarrow\infty}~\pi_{\epsilon}
$. The uniqueness property of the invariant measure and the ergodicity properties of the process $(A^{\epsilon}_t)_{t\geq 0}$ require more sophisticated stochastic analysis.

Assume that the process $A^{\epsilon}$ is mean ergodic, in the sense that
\begin{equation}\label{ref-ergodic-time}
\frac{1}{t}\int_0^{t}~\Vert A^{\epsilon}_{s}-\widehat{A}^{\epsilon}_{\infty}\Vert~ds\longrightarrow_{t\rightarrow\infty} 0\quad\mbox{\rm a.s.}\quad\mbox{\rm with}\quad\widehat{A}^{\epsilon}_{\infty}:=\int_{\Ma_r}~\Lambda~\pi_{\epsilon}(d\Lambda)
\end{equation}
By the convexity property of the log-norm we have the almost sure convergence result
\begin{equation}\label{ref-ergodic-time-2}
\frac{1}{t}\log{\Vert \Ea_{t}(A^{\epsilon})\Vert}\leq \frac{1}{t}\int_0^{t}~\mu(A^{\epsilon}_{s})~ds~~\longrightarrow_{t\rightarrow\infty}  ~~\mu(\widehat{A}^{\epsilon}_{\infty})
\leq \int_{\Ma_r}~\mu(\Lambda)~\pi_{\epsilon}(d\Lambda)
\end{equation}
In this situation, the condition $\mu(\widehat{A}^{\epsilon}_{\infty})<0$ for some $\epsilon\in ]0,1[$ ensures that for any $\nu\in ]0,1[$ there exists a random time $\tau^{\epsilon}_{\nu}$ such that 
\begin{equation}\label{ref-ergodic-time-3}
t \,\geq\, \tau^{\epsilon}_{\nu} ~~\quad\Longrightarrow\quad~~ \frac{1}{t}\log{\Vert \Ea_{t}(A^{\epsilon})\Vert}\,\leq\, (1-\nu)\,\mu(\widehat{A}^{\epsilon}_{\infty}) \,<\, 0
\end{equation}
This provides a natural condition under which the semigroup $\Ea_{t}(A^{\epsilon})$ is almost surely exponentially stable in terms of the long time behaviour of the stochastic process $A^{\epsilon}$. This approach is related to so-called averaging methods for random dynamical systems of the form (\ref{stochastic-OU-linear}); e.g. see \cite{Blankenship77,Geman79,solo69}. 

For instance, in the context of the ensemble Kalman-Bucy filter, the stochastic Riccati diffusion (\ref{ric-intro}) forgets exponentially fast its initial condition (cf. Theorem~2.4 in~\cite{2018arXiv180800235}), and if $\mu({A}_{\infty})<0$ holds under strong enough controllability and observability conditions, then the log-norm $\mu(\widehat{A}^{\epsilon}_{\infty})$ can be made as close as desired to $\mu({A}_{\infty})$ with a sufficiently small fluctuation parameter $\epsilon$. It follows that the random matrices $A^{\epsilon}_t$ in (\ref{KB-intro}) may satisfy (\ref{ref-ergodic-time}) and (\ref{ref-ergodic-time-2}) with $\mu(\widehat{A}^{\epsilon}_{\infty})<0$; but this requires a (difficult to quantify) sufficiently large number of samples (or $\epsilon$ sufficiently small). Also, unfortunately, the decays rate to equilibrium of the Riccati diffusion (\ref{ric-intro}) are difficult to quantify, and they depend on the fluctuation parameter $\epsilon$. Therefore, even if this approach could be applied to infer the random semigroup $\Ea_{t}(A^{\epsilon})$ associated with $A^{\epsilon}_t$ in (\ref{KB-intro}) is almost surely exponentially stable, this analysis doesn't provide any useful information on the decay rates of the semigroup.

Unfortunately, in general applications, the mean ergodic property (\ref{ref-ergodic-time}) of the random matrix valued process $A^{\epsilon}$, as well as the condition $\mu(\widehat{A}^{\epsilon}_{\infty})<0$, can be difficult to check. It is also generally difficult to quantify with any precision the almost sure convergence rate in the limiting equation (\ref{ref-ergodic-time-2}). Moreover, the rather crude asymptotic analysis discussed above doesn't provide any useful information on the relaxation time $\tau^{\epsilon}_{\nu}$ introduced (\ref{ref-ergodic-time-3}). Finally, without additional information this strategy also cannot be used to estimate the $n$-th moments of the process $X^{\epsilon,0}_t$ as in (\ref{ref-intro-perfect-est-conseq}).

In another direction, assume now that $q_1\mathrm{Id}\leq Q_t\leq q_2\mathrm{Id}$ is a possibly random solution of the Lyapunov equation
$$
\partial_tQ_t+A^{\epsilon}_tQ_t+Q_t(A^{\epsilon}_t)^{\prime} \,\leq\, -q_3\,\mathrm{Id}
$$
for some $q_1,q_2,q_3>0$. Then, we have the almost sure contraction estimate
\begin{align*}
\partial_t \left[\Ea_t(A^{\epsilon})^{\prime}\,Q_t\,\Ea_t(A^{\epsilon})\right] &\,\leq \, -(q_3/q_2)\,\,\Ea_t(A^{\epsilon})^{\prime}\,Q_t\,\Ea_t(A^{\epsilon})\\
 & \Longrightarrow\quad \Ea_t(A^{\epsilon})^{\prime}\Ea_t(A^{\epsilon}) \,\leq \, (q_2/q_1)\,e^{-(q_3/q_2)t}~\mathrm{Id}\\
 & \Longrightarrow\quad  \Vert \Ea_t(A^{\epsilon})\Vert \,\leq \, 
\alpha\,e^{-\beta t}~~ \mbox{\rm with parameters}~~(\alpha,\beta)=(\sqrt{q_2/q_1},~2^{-1} q_3/q_2)
\end{align*}
Conversely, suppose that $\Vert A^{\epsilon}_t\Vert\leq \gamma$ is almost surely uniformly bounded and the almost sure contraction estimate just given is satisfied. In this case, for any matrix valued process $R_t$ satisfying
$$
	r_1\mathrm{Id}\leq R_t\leq r_2\mathrm{Id}
$$ 
we have 
$$
Q_t \,:=\, \int_t^{\infty}\Ea_{t,u}(A^{\epsilon})^{\prime}R_u~\Ea_{t,u}(A^{\epsilon})\,du \,\leq\, r_2\alpha/\beta
\quad\Longrightarrow\quad \partial_tQ_t+(A^{\epsilon}_t)^{\prime}Q_t+Q_tA_t^{\epsilon}\,=\,-R_t \,\leq\, -r_1\,\,
$$
as well as
$$
\Ea_{t,s}(A^{\epsilon})^{\prime}\,\Ea_{t,s}(A^{\epsilon}) ~\geq~ e^{-\int_s^t\lambda_{max}(-(A^{\epsilon})_{sym})}~\geq~ e^{-\gamma (t-s)}\quad \Longrightarrow\quad
Q_t~\geq~ (r_1/\gamma)\, \mathrm{Id}
$$
This provides sufficient and necessary conditions under which  $\Ea_{t}(A^{\epsilon})$ is almost surely exponentially stable in terms of the existence of Lyapunov functions. 

Unfortunately, the design of Lyapunov functions for nonlinear stochastic diffusions in matrix spaces is a difficult task. Stability (in the mean) of linear random dynamical equations of the form (\ref{stochastic-OU-linear}) via Lyapunov methods was considered in early work by Bertram and Sarachik \cite{Bertram1959} and in \cite{kats60}. However, application of this method in the mean is also typically not practical \cite{Bertram1959,Kozin69}. 

For example, in view of the complexity of the stochastic Riccati equation in (\ref{ric-intro}), the Lyapunov technique discussed above doesn't apply to the random matrices $A^{\epsilon}_t$ in (\ref{KB-intro}).

 The aim of this work is to provide some quantitative stability properties under weaker conditions.
Namely, we shall consider the following regularity conditions:
\begin{hypothesis}[$H_1$]
Suppose $H_0$ holds and,
$$
\forall n\geq 1~~~~\exists\epsilon_n\in [0,1]~~~~ \mbox{\rm such that}~~~~
\forall \epsilon\in [0,\epsilon_n] ~~~~\mbox{\rm  we have~~~~}
\sup_{t\geq 0}\,{\EE\left(\Vert A_t-A^{\epsilon}_t\Vert^n\right)^{1/n}}\leq c_n\,\epsilon
$$
We also define throughout $\epsilon_n(\nu):=\epsilon_n\wedge \left[\nu^{1/n}c_0/(4c_n)\right]$ indexed by $n\geq 1$ and $\nu\in]0,1[$.
\end{hypothesis}

\begin{hypothesis}[$H_2$]
Suppose $H_0$ holds and,
\begin{equation*}
\forall n\geq 1 ~~~~\mbox{\rm and}~~~~ \forall \epsilon\in [0,1]~~~~\mbox{\rm we have the uniform estimates}
\end{equation*}
\begin{equation*}
 \displaystyle\sup_{t\geq 0}\,{\EE\left(\Vert A_t-A_t^{\epsilon}\Vert^n\right)^{1/n}} \,\leq\, c_n\,\epsilon~~~~\mbox{\rm and}~~~~
\sup_{t\geq 0}\,{\EE\left(\Vert A^{\epsilon}_t\Vert^n\right)^{1/n}}\,\leq\, d_1+\epsilon\,d_2\,n^{1/2}
\end{equation*}
\end{hypothesis}

In the above hypotheses, $d_1\leq d_2$, $c_n$ correspond to any non-decreasing collection of finite non-negative constants. Observe that
$(H_2)\Longrightarrow (H_1)$. These hypotheses both correspond essentially to some level of fluctuation control on $A^{\epsilon}_t$, around a (time-varying) deterministic limit that is itself stabilising. This is a rather natural condition to propose, and quite weak. Indeed, in the context of mean-ergodic and mean-stable systems, e.g. as discussed in the literature, see (\ref{ref-ergodic-time}), we can immediately view these hypotheses as a significant relaxation. In this context, we only require some fluctuation control on the matrices themselves at each time and at any moment these matrices may themselves be unstable.

\begin{hypothesisRemark}[$H_0''$]
Conditions $(H_1)$ and $(H_2)$ may also be thought of as weak conditions on the moment continuity between $A_t$ and $A^{\epsilon}_t$ in terms of $\epsilon$. It may be possible to relax the log-norm condition in $(H_0)$ whenever we have a stability estimate (e.g. like (\ref{stabilityestimate1})) that respects these continuity properties. Unfortunately, apart from very particular cases (e.g. time-invariant flows, or time-varying commuting matrices) we are aware of no general relaxation to the log-norm condition in the time-varying setting. 
\end{hypothesisRemark}

For a deterministic flow of perturbed matrices $A^{\epsilon}_t$, these two conditions coincide when the fluctuation parameter is chosen sufficiently small. In this case, the deterministically perturbed semigroup $\Ea_{s,t}(A^{\epsilon})$ is exponentially stable for small perturbations.

Returning to our running illustrative example, it follows by Theorem 2.3 in~\cite{2018arXiv180800235}, the stochastic matrices (\ref{KB-intro})  associated with the various classes of ensemble Kalman-Bucy filter mentioned in the introduction satisfy condition $(H_1)$. In addition, by Theorem~ 2.3 in~\cite{2018arXiv180800235}, condition $(H_2)$ is met when the matrix $\Va$ introduced in (\ref{sigma-intro}) is null; which corresponds to a particular form of ensemble filter. These conditions hold in this example for any reasonable choice of sample size (i.e. any $\epsilon$).

We end this section with some comments on the above regularity conditions. Firstly, we emphasize that these conditions don't require any ergodic property on the process $A_t^{\epsilon}$, nor any conditions on the limiting log-norm $\mu(\cdot)$. 

When $(H_1)$ is satisfied we have the estimate
$$
\sup_{t\geq 0}\,{\EE\left(\Vert A^{\epsilon}_t\Vert^n\right)^{1/n}}~\leq~  (a+\Vert A_{\infty}\Vert)+c_n\,\epsilon
$$
This ensures that the l.h.s. condition in (\ref{hyp-unif-moments}) is met with $\epsilon_{0,n}=\epsilon_n$ when $(H_1)$ is met. When $(H_2)$ is satisfied then the l.h.s. condition in (\ref{hyp-unif-moments}) is trivially met with $\epsilon_{0,n}=1$. 

This also shows that the r.h.s. condition in $(H_2)$ is met as soon as the l.h.s. condition in $(H_2)$ is satisfied with $c_n=n^{1/2}$.
The latter is often difficult to check for nonlinear diffusion approximation models since the fluctuation analysis of the $n$-th error moments often combine Burkh\"older-Davis-Gundy-type inequalities involving a square root parameter $n^{1/2}$, with the estimation of $n$-th order type moments of $A^{\epsilon}_t$.

When $B^{\epsilon}$ depends on $A^{\epsilon}$ w.r.t. some polynomial type function, the r.h.s. condition in (\ref{hyp-unif-moments}) can be readily checked using the moments estimates on $A^{\epsilon}$ just discussed.

Finally, we may illustrate the satisfaction of $(H_2)$ with some additional examples: 
\begin{itemize}
\item Firstly notice that $(H_2)$ is satisfied for the spectral norm and sub-Gaussian fluctuations. 

To be more precise, we set
\begin{equation*}
(\Delta\mathrm{A})^{\epsilon}_t:=\epsilon^{-1}\left[A^{\epsilon}_t-A_t\right]~~\Longleftrightarrow ~~A^{\epsilon}_t=A_t+\epsilon\,(\Delta\mathrm{A})^{\epsilon}_t
\end{equation*}
In this notation, we have
$$
\sup_{t\geq 0}{\EE\left(\Vert (\Delta\mathrm{A})^{\epsilon}_t\Vert^n\right)^{1/n}}\leq d_2\,n^{1/2}~~\Longrightarrow~~
(H_2)\quad \mbox{\rm with}\quad c_n=d_2\,n^{1/2}\quad \mbox{\rm and}\quad d_1=a+\Vert A_{\infty}\Vert
$$
Also observe that $(H_2)$ is met for fluctuation matrices with sub-Gaussian entries; that is, when the following condition is met
$$
\forall 1\leq i,j\leq r,\quad \sup_{t\geq 0}\,{\EE\left(\left[(\Delta\mathrm{A})^{\epsilon}_t(i,j)\right]^n\right)^{1/n}}\,\leq\, d_2\,n^{1/2}~~~\Longrightarrow~~~ (H_2)
$$
\item Let $Z^{\epsilon}_t$ be some collection of stochastic processes indexed by $\epsilon\in [0,1]$, independent of $W_t$ and taking values in  some Banach space $\mathcal{B}$, equipped with some norm $\Vert \cdot\Vert_{\mathcal{B}}$.

We further assume that $Z_t:=Z^{0}_t$ is a deterministic process that converges exponentially fast to some value $Z_{\infty}$ when $t\rightarrow\infty$.

Suppose that the pair of matrices
$(A^{\epsilon}_t,B^{\epsilon}_t)$ are given by
$$
A^{\epsilon}_t=C\left(Z^{\epsilon}_t\right)\quad \mbox{\rm and}\quad
(B^{\epsilon}_t)^{1/2}=D\left(Z^{\epsilon}_t\right)\quad\Longrightarrow\quad dX^{\epsilon}_t=C\left(Z^{\epsilon}_t\right)\,X^{\epsilon}_t\,dt+D\left(Z^{\epsilon}_t\right)\,dW_t
$$
for some Lipschitz drift functions $C\,:\,\mathcal{B}\mapsto \Ma_r$ and some positive diffusion map $D\,:\,\mathcal{B}\mapsto  \Sa^0_r$ with polynomial growth.

 In this case, condition (\ref{hyp-unif-moments}) is satisfied as soon as the process $Z^{\epsilon}_t$ has uniformly bounded moments at any order w.r.t. the time parameter. Condition $H_0$ is also met with
$$
A_t=C\left(Z_t\right)\quad \mbox{\rm and}\quad A_{\infty}=C\left(Z_{\infty}\right)\quad \mbox{\rm as soon as}\quad\mu\left(C\left(Z_{\infty}\right)\right)<0
$$
In addition, condition $(H_1)$, resp.   $(H_2)$, is satisfied as soon as the processes $(Z^{\epsilon}_t,Z_t)$ satisfy the regularity condition defined as 
$(H_1)$, resp.  $(H_2)$, by replacing $(A^{\epsilon}_t,A_t)$ by $(Z^{\epsilon}_t,Z_t)$ and the matrix norm by the norm $\Vert \cdot\Vert_{\mathcal{B}}$.
\end{itemize}

\subsection{A Stochastic Stability Theorem}

In this section we state the main result of this work along with a number of ancillary corollaries of interest on their own.

Our main result takes the following form.
\begin{theo}\label{theo-A-epsilon}~

$\bullet~$Suppose the fluctuation estimates in $(H_1)$ are satisfied for $n=2$. Then, for any time horizon 
$s\geq 0$, and any parameter $\nu\in]0,1[$ we have
 \begin{align}
 t\geq \frac{2}{\nu}\frac{a}{bc_0}\quad &\mbox{and} \quad \epsilon\leq \epsilon_2\wedge \left(\frac{\nu}{2}
\frac{c_0}{c_2}\right) \nonumber \\
& \Longrightarrow\quad
\log{\Vert \Ea_{s,s+t}(\overline{A}^{\,\epsilon})\Vert}\,\vee\, \EE\left[\log{\Vert \Ea_{s,s+t}({A}^{\,\epsilon})\Vert}\right] \,\leq\, (1-\nu)\,\mu(A_{\infty})\,t \label{EA-epsilon-H1-}
\end{align} 
where $\overline{A}^{\,\epsilon}:t\in\RR_+\mapsto \overline{A}^{\,\epsilon}_t:=\EE(A^{\epsilon}_t)$ is the averaged process. 

$\bullet~$Assume $(H_1)$ is satisfied. Then, for any increasing sequence of times $0\leq s \leq t_k\uparrow_{k\rightarrow\infty}\infty$, the probability of the following event
 \begin{equation}\label{EA-epsilon-H1}
\limsup_{k\rightarrow\infty}\frac{1}{t_k}\log{\Vert \Ea_{s,t_k}(A^{\epsilon})\Vert} \,< \, \frac{1}{2}\,\mu(A_{\infty})\quad \mbox{is greater than $1-\nu$}
 \end{equation}
for any $\nu\in]0,1[$, as soon as  $\epsilon$ is chosen so that $\epsilon \leq \epsilon_n(\nu)$ for some $n\geq 1$.

$\bullet~$Now suppose hypothesis $(H_2)$ is satisfied. Then, for any $n\geq 1$, any fluctuation parameter $\epsilon\leq \epsilon_{2,n}$, and any time horizon $s\geq 0$ we have
 \begin{equation}\label{EA-epsilon-H2}
T_n\leq  t\leq T_{n}^{\epsilon}~~\Longrightarrow~~\frac{1}{t}\log{\EE\left(\Vert \Ea_{s,s+t}(A^{\epsilon})\Vert^{n}\right)}\,\leq \,
 \frac{n}{4}\,\mu(A_{\infty})
 \end{equation}
 where 
 \begin{equation}\label{def-Tn-eps}
T_n:=\frac{4}{c_0}~\log{\left(1+\frac{c_{2n}}{c_0}~2^{2+\frac{1}2n}\right)}\quad \mbox{\rm and}\quad T_{n}^{\epsilon}:=
\frac{\log{(1/\epsilon^2)}}{2((4e+(2e)^{1/2})~(d_1\vee d_2))+c_0}\wedge \frac{1/\epsilon^2}{2d_2n}
 \end{equation}
and $\epsilon_{2,n}$ is the largest parameter such that $T_n<T_{n}^{\epsilon}$. Note that $T_{n}^{\epsilon}\longrightarrow_{\epsilon\rightarrow 0}\infty$. 
\end{theo}

The proof of the above theorem is provided in Section~\ref{theo-A-epsilon-proof}. We illustrate the impact of the above theorem with a series of corollaries outlined in the subsequent subsection. 

\subsection{Corollary Results} 

Firstly, we consider a collection of corollaries under the hypothesis $(H_1)$. The first corollary is a consequence of the Borel-Cantelli's lemma applied to (\ref{EA-epsilon-H1}) in Theorem \ref{theo-A-epsilon}.

\begin{cor}\label{corBC-H1}
Assume $(H_1)$ is satisfied. Then, for any $s\geq 0$, any increasing sequence of time horizons $t_{k_1}\uparrow_{{k_1}\rightarrow\infty}\infty$ and any sequence $\varepsilon_{k_2}\downarrow_{{k_2}\rightarrow\infty} 0$ such that $\sum_{{k_2}\geq 1}\varepsilon_{k_1}^{\,n}<\infty$ for some $n\geq 1$, we have the almost sure Lyapunov estimate
 \begin{equation}\label{EA-epsilon-H1-dd}
\limsup_{{k_2}\rightarrow\infty}\limsup_{{k_1}\rightarrow\infty}\frac{1}{t_{k_1}}\,\log{\Vert \Ea_{s,s+t_{k_1}}(A^{\varepsilon_{k_2}})\Vert}\,<\, \frac{1}{2}\,\mu(A_{\infty})
 \end{equation}
\end{cor}

The next two results provide some reformulation of the supremum limit estimates stated in (\ref{EA-epsilon-H1}) and (\ref{EA-epsilon-H1-dd}) in terms of random relaxation time horizons and random relaxation-type fluctuation parameters.

The first of these two results shows that with a high probability, the semigroup $\Ea_{s,t}(A^{\epsilon})$ is stable after some 
possibly random relaxation time horizon, as soon as $\epsilon$ is chosen sufficiently small.

\begin{cor}\label{corH1ref}
Assume $(H_1)$ holds. Then, for any increasing sequence of times $0\leq s \leq t_k\uparrow_{k\rightarrow\infty}\infty$, the probability of the following event,   
 \begin{align}
 \left\{\begin{array}{l}
 \forall 0<\nu_2\leq 1~~~ \exists l\geq 1 ~~~\mbox{such~that}~~~ \forall k\geq l~~~\mbox{it~holds~that~} \\~
 \\
 \qquad\qquad\qquad\qquad\qquad\qquad\qquad\qquad\displaystyle\frac{1}{t_k}\log{\Vert \Ea_{s,s+t_k}(A^{\epsilon})\Vert} \,\leq\,   \frac{1}{2}\,(1-\nu_2)\,\mu(A_{\infty})
 \end{array}\right\} \label{EA-epsilon-H1-cor}
\end{align}
 is greater than $1-\nu_1$, for any  $\nu_1\in]0,1[$, as soon as  $\epsilon$ is chosen so that $\epsilon \leq \epsilon_n(\nu_1)$ for some $n\geq 1$.
\end{cor}

The next result takes this one step further for an almost sure result that comes into effect after some random time and with some sufficiently small fluctuation parameter (where sufficiency in this case is also random). 

\begin{cor}\label{corH1ref-v2}
Assume $(H_1)$ is satisfied. Consider any $s\geq 0$, any increasing sequence of time horizons $t_k\uparrow_{{k_1}\rightarrow\infty}\infty$, and any sequence $\varepsilon_{k_2}\downarrow_{{k_2}\rightarrow\infty} 0$ such that $\sum_{{k_2}\geq 1}\varepsilon_{k_2}^{\,n}<\infty$ for some $n\geq 1$. Then, we have the almost sure Lyapunov estimate, 
\begin{align}
 \left\{\begin{array}{l}
 \forall 0<\nu\leq 1~~~ \exists l_1,l_2\geq 1 ~~~\mbox{such~that}~~~ \forall k_1\geq l_1,~\forall k_2\geq l_2~~~\mbox{it~holds~that~} \\~
 \\
 \qquad\qquad\qquad\qquad\qquad\qquad\qquad\qquad\displaystyle \frac{1}{t_{k_1}}\log{\Vert \Ea_{s,s+t_{k_1}}(A^{\varepsilon_{k_2}})\Vert} \,\leq\,   \frac{1}{2}\,(1-\nu)\,\mu(A_{\infty})
\end{array}\right\} \label{EA-epsilon-H1-cor-v2}
\end{align}

\end{cor}

The formulation in the preceding corollary underlines the fact that after some random time (i.e. defined in terms of $l_1$), and given some randomly sufficiently small perturbation (defined in terms of $l_2$) the system (\ref{stochastic-OU}) will almost surely achieve asymptotic (exponential) stability. We have no control over the parameters $l_1$ and $l_2$ which depend on the randomness in a realisation of $A^{\varepsilon_{k_2}}_t$.

The next result concerns the stability of the process (\ref{stochastic-OU}) itself. 

\begin{cor}\label{lem-tech-Aat-1}
Assume $(H_1)$ holds. Then, for any increasing sequence of time horizons $t_k\uparrow_{k\rightarrow\infty}\infty$ and any $x_1\not=x_2$, the probability of the following event
 \begin{equation}\label{EA-epsilon-H1-cor-v3}
\limsup_{k\rightarrow\infty}\frac{1}{t_k}\log{\Vert X_{t_k}^{\epsilon,x_1}-X_{t_k}^{\epsilon,x_2}\Vert} \,<\, \frac{1}{2}\,\mu(A_{\infty})\quad \mbox{is greater than $1-\nu$}
 \end{equation}
for any  $\nu\in]0,1[$, as soon as  $\epsilon$ is chosen so that $\epsilon \leq  \epsilon_n(\nu)$ for some $n\geq 1$.
\end{cor}
The preceding corollary is a direct consequence of the decomposition
$$
X_t^{\epsilon,x_1}-X_t^{\epsilon,x_2}=\Ea_{t}(A^{\epsilon})~(x_1-x_2)~~\Longrightarrow~~\Vert X_t^{\epsilon,x_1}-X_t^{\epsilon,x_2}\Vert\leq
\Vert \Ea_{t}(A^{\epsilon})\Vert~\Vert x_1-x_2\Vert
$$
Note that (\ref{EA-epsilon-H1-cor-v3}) in Corollary \ref{lem-tech-Aat-1} is analogous to (\ref{EA-epsilon-H1}) in Theorem \ref{theo-A-epsilon} but at the level of the process (\ref{stochastic-OU}) itself. Analogous results to Corollaries \ref{corBC-H1}, \ref{corH1ref}, and \ref{corH1ref-v2} at the level of the process (\ref{stochastic-OU}) follow immediately.

Next, we consider a collection of corollaries under the stronger hypothesis $(H_2)$. Firstly, given $(H_2)$, we highlight a fact immediate from (\ref{EA-epsilon-H2}) and (\ref{def-Tn-eps}), that for any $n\geq 1$, any $s\geq0$, we have
$$
\limsup_{\epsilon\rightarrow0}\,
\frac{1}{T_{n}^{\epsilon}}\,\log{\EE\left(\Vert \Ea_{s,s+T_{n}^{\epsilon}}(A^{\epsilon})\Vert^{n}\right)} \,\leq \,
 \frac{n}{4}\,\mu(A_{\infty})
$$

The next result provides a fluctuation-type analysis. 

\begin{cor}
Suppose $(H_2)$ is satisfied. For any $\epsilon\leq \epsilon_{2,2n}$ we have the fluctuation estimate
\begin{equation}
T_{2n}\leq  s\leq t\leq T_{2n}^{\epsilon} \quad \Longrightarrow \quad
\epsilon^{-1}\,\EE\left(\Vert \Ea_{s,t}(A^{\epsilon})-\Ea_{s,t}(A)\Vert^n\right)^{1/n} \,\leq\, c_n+4e^{a/b}\,c_{2n}/c_0
\end{equation}
\end{cor}
\proof
The first assertion follows immediately from (\ref{EA-epsilon-H2}). To check the second assertion, we use the decomposition
\begin{align*}
& \partial_t \left[\Ea_{s,t}(A^{\epsilon})-\Ea_{s,t}(A)\right]\,=\,A_t\left[\Ea_{s,t}(A^{\epsilon})-\Ea_{s,t}(A)\right]+(A^{\epsilon}_t-A_t)\,\Ea_{s,t}(A^{\epsilon}) \\ 
&  \quad \Longrightarrow \quad \Ea_{s,t}(A^{\epsilon})-\Ea_{s,t}(A) \,=\, (A^{\epsilon}_s-A_s)+\int_s^t\Ea_{u,t}(A)\,(A^{\epsilon}_u-A_u)\,\Ea_{s,u}(A^{\epsilon})\,du\\ 
& \quad \Longrightarrow \quad \epsilon^{-1}\,\EE\left(\Vert \Ea_{s,t}(A^{\epsilon})-\Ea_{s,t}(A)\Vert^n\right)^{1/n} \,\leq\, c_n+c_{2n}e^{a/b}
\int_s^t\,e^{(t-u)\mu(A_{\infty})}\,\EE\left(\Vert\Ea_{s,u}(A^{\epsilon})\Vert^{2n}\right)^{1/(2n)}\,du
\end{align*}
This implies that
\begin{align*} 
\quad &  T_{2n}\leq  s\leq t\leq  T_{2n}^{\epsilon}  \\
 & \quad\quad~~ \Longrightarrow  \quad  
\epsilon^{-1}\,\EE\left(\Vert \Ea_{s,t}(A^{\epsilon})-\Ea_{s,t}(A)\Vert^n\right)^{1/n}\,\leq \,c_n+\frac{4c_{2n}}{c_0}e^{a/b}\,e^{(t-s)\mu(A_{\infty})/4}\,\left(1-e^{3(t-s)\mu(A_{\infty})/4}\right)
\end{align*}
The proof of the corollary is complete.\qed

The next corollary concerns stability in the mean, at the level of the process (\ref{stochastic-OU}) itself, and guaranteed over a relevant (computable) deterministic time interval.  

The next corollary establishes and makes precise the relationship alluded to in prior discussion, i.e. (\ref{ref-intro-perfect-est})~$\Rightarrow$~(\ref{ref-intro-perfect-est-conseq}). It is based on the fact that under $(H_2)$ alone, the result (\ref{EA-epsilon-H2}) in Theorem \ref{theo-A-epsilon} establishes an estimate of the form (\ref{ref-intro-perfect-est}), at least over an interval (which can be chosen as large as needed by reducing the fluctuation parameter).

\begin{cor}\label{lem-tech-Aat}
Suppose $(H_2)$ holds. Then, for any $n\geq 1$, $\epsilon\leq \epsilon_{2,n}$ and any time horizon $t$ such that $T_n\leq t\leq T^{\epsilon}_n$, we have the contraction inequality,
\begin{equation}\label{uniformly-Lipschitz-moments-OU-Lip}
{\EE\left(\Vert X_t^{\epsilon,x_1}-X_t^{\epsilon,x_2}\Vert^n\right)}^{1/n} \,\leq\, \exp{\left[\mu(A_{\infty})t/4\right]}~\Vert x_1-x_2\Vert~
\end{equation}

In addition, for any $n\geq 2$ and any $\epsilon\leq \epsilon_{2,n}\wedge \epsilon_{1,2n}$, we have the moment estimates
\begin{equation}\label{uniformly-Lipschitz-moments-OU}
T_n\vee T_{2n} \leq t\leq  T^{\epsilon}_{2n}\quad\Longrightarrow\quad\EE\left(\Vert X_t^{\epsilon,x}\Vert^n\right)^{1/n} \,\leq\,
 \exp{\left[\mu(A_{\infty})t/4\right]}\,\Vert x\Vert+ \kappa_n
\end{equation}
for some finite constant $\kappa_n$ whose value only depends on the parameter $n$ (and possibly on $r$).
\end{cor}
\proof
Observe that
$$
\Vert X_t^{\epsilon,x}-X_t^{\epsilon,y}\Vert \,\leq\, \Vert \Ea_{t}(A^{\epsilon})\Vert\,\Vert x-y\Vert\quad\mbox{\rm and}\quad
\Vert X_t^{\epsilon,x}\Vert \,\leq\, \Vert X_t^{\epsilon,0}\Vert+\Vert \Ea_{t}(A^{\epsilon})\Vert\,\Vert x\Vert
$$
The estimate (\ref{uniformly-Lipschitz-moments-OU-Lip}) is a direct consequence of (\ref{EA-epsilon-H2}) and the l.h.s. estimate in the above display. On the other hand, we have
\begin{align*} 
\Vert X_t^{\epsilon,0}\Vert^{2n}&~=~\Vert \int_0^t~\Ea_{s,t}(A^{\epsilon})~(B^{\epsilon}_s)^{1/2}~dW_s\Vert^{2n}~=~\left[\sum_{1\leq i\leq r}
\left[\sum_{1\leq j\leq r}\int_0^t~\left[\Ea_{s,t}(A^{\epsilon})~(B^{\epsilon}_s)^{1/2}\right]_{i,j}~dW_s^j\right]^2\right]^n\\
&~ \leq ~r^{3n-2}\,\sum_{1\leq i,j\leq r}
\left[\int_0^t~\left[\Ea_{s,t}(A^{\epsilon})~(B^{\epsilon}_s)^{1/2}\right]_{i,j}~dW_s^j\right]^{2n}
\end{align*}
By the Burkh\"older-Davis-Gundy inequality (e.g. Proposition 4.2 in~\cite{barlow1982semi}), for any $n\geq 1$ we have
\begin{eqnarray*}
\EE\left[\Vert  X_t^{\epsilon,0}\Vert^{2n}~\vert~\Fa^{\epsilon}_t\right]&\leq&  c\,r^{3n-2}\,(2n)^{n}\,\sum_{1\leq i,j\leq r}
\left[\int_0^t~\left[\Ea_{s,t}(A^{\epsilon})~(B^{\epsilon}_s)^{1/2}\right]_{i,j}^2~ds\right]^{n}
\end{eqnarray*}
for some universal constant $c$. This yields
\begin{eqnarray*}
\EE\left[\Vert  X_t^{\epsilon,0}\Vert^{2n}~\vert~\Fa^{\epsilon}_t\right]
&\leq& c\,(2n)^{n}\,\left[\int_0^t~\tr\left[\Ea_{s,t}(A^{\epsilon})~B^{\epsilon}_s~\Ea_{s,t}(A^{\epsilon})^{\prime}\right]~ds\right]^{n}\\
&\leq & c\,(2n)^{n}\,\left[\int_0^t~\Vert \Ea_{s,t}(A^{\epsilon})\Vert^2~\tr\left[B^{\epsilon}_s\right]~ds\right]^{n}
\end{eqnarray*}
where $c$ may vary from line to line (and depend on $r$ but not on $n$). Combining (\ref{hyp-unif-moments}) with the generalized Minkowski inequality and Cauchy-Schwartz inequality we check the estimate
\begin{align*} 
\EE\left[\Vert X_t^{\epsilon,0}\Vert^{2n}\right]^{1/n}
\,\leq\,  c^{1/n}\,n\,\rho_{2n}~\int_0^t~\EE\left(\Vert \Ea_{s,t}(A^{\epsilon})\Vert^{4n}\right)^{1/(2n)}\,ds
\end{align*}
Assume that condition $(H_2)$ is satisfied. Observe that for any $T_{4n}\leq t\leq T_{4n}^{\epsilon}$ we have
\begin{align*}
&\int_0^{t}~\EE\left[\Vert \Ea_{s,t}(A^{\epsilon})\Vert^{4n}\right]^{1/(2n)}\,ds\\
& \qquad~~ =~
\int_0^{T_{4n}}~\EE\left[\Vert \Ea_{(t-s),(t-s)+s}(A^{\epsilon})\Vert^{4n}\right]^{1/(2n)}\,ds+
\int_{T_{4n}}^{t}~\EE\left[\Vert \Ea_{(t-s),(t-s)+s}(A^{\epsilon})\Vert^{4n}\right]^{1/(2n)}\,ds
\end{align*}
By Theorem~\ref{theo-A-epsilon}, we have
\begin{eqnarray*}
\int_{T_{4n}}^{t}\,\EE\left[\Vert \Ea_{(t-s),(t-s)+s}(A^{\epsilon})\Vert^{4n}\right]^{1/(2n)}\,ds&\leq&\int_{T_{4n}}^{T_{4n}^{\epsilon}}~\exp{\left[
\frac{s}{2}\,\mu(A_{\infty})\right]}\,ds\\
 &=&\frac{2}{\vert\mu(A_{\infty})\vert}\,\left[\exp{\left[\frac{T_{4n}}{2}\,\mu(A_{\infty})\right]}-\exp{\left[\frac{T_{4n}^{\epsilon}}{2}\,\mu(A_{\infty})\right]}\right]
\end{eqnarray*}
In this situation, using (\ref{ref-nE-exp-Aa-1}) we also have
\begin{align*}
& \int_0^{T_{4n}}\,\EE\left[\Vert \Ea_{(t-s),(t-s)+s}(A^{\epsilon})\Vert^{4n}\right]^{1/(2n)}\,ds\\
& \qquad~~\leq~ \int_0^{T_{4n}}\,\left[ \frac{1}{2}\,e^{8d_1ns}+\frac{e}{2}\,\sqrt{\frac{e}{\pi}}
\left(4(1+(2e)^{1/2}d_2ns\epsilon)\,e^{(8e^{1/2}d_2ns\epsilon)^2}-1\right)\right]^{1/(2n)}\,ds\\
& \qquad~~\leq~  \frac{1}{4d_1}\,\left[e^{4d_1T_{4n}}-1\right]+\left[2e\,\sqrt{\frac{e}{\pi}}
(1+(2e)^{1/2}d_2nT_{4n}\epsilon)\right]^{1/(2n)}\,\frac{1}{e4d_2}\left[e^{e4d_2T_{4n}}-1\right]
\end{align*}
 In the last assertion we have used the fact that
 $
4nd_2T_{4n}\epsilon^2\leq 1/2
 $.
These estimates imply that
$$
\sup_{0\leq t\leq T_{4n}^{\epsilon}}{\EE\left[\Vert X_t^{\epsilon,0}\Vert^{2n}\right]^{1/(2n)}}
\displaystyle\leq  c\,n^{1/2}\,(\rho_{2n}\delta^{\epsilon}_{2n})^{1/2}
$$
with
$$
\delta_{2n}^{\epsilon} \,:=\, \frac{1}{4d_1}\,\left[e^{4d_1T_{4n}}-1\right]+\left[
1+d_2nT_{4n}\epsilon\right]^{1/(2n)}\,\frac{1}{d_2}\left[e^{4ed_2T_{4n}}-1\right]
$$
Using (\ref{EA-epsilon-H2}) we conclude that 
$$
T_n\vee T_{2n}\leq t\leq T^{\epsilon}_{2n}\quad\Longrightarrow\quad
\EE\left(\Vert X_t^{\epsilon,x}\Vert^n\right)^{1/n}\leq c\,n^{1/2}\,(\rho_{n}\delta^{\epsilon}_{n})^{1/2}+
e^{\mu(A_{\infty})t/4} \,\Vert x\Vert
$$
for any $n\geq 2$. This ends the proof of the corollary. \qed

\section{Proof of Theorem~\ref{theo-A-epsilon}}\label{theo-A-epsilon-proof}
The proof of Theorem~\ref{theo-A-epsilon} is based on the following technical lemma.
\begin{lem}
Assume that the r.h.s. estimates in $(H_2)$ are satisfied. Then, for any $n\geq 1$, $\epsilon\in [0,1]$ and $s,t\geq 0$ we have
the estimate
\begin{equation}\label{ref-nE-exp-Aa-1}
\EE\left(\Vert  \Ea_{s,s+t}(A^{\epsilon})\Vert^n\right) \,\leq\,   \frac{1}{2}\,e^{2d_1nt}+\frac{e}{2}\,\sqrt{\frac{e}{\pi}}
\left((1+(2e)^{1/2}\,d_2nt\epsilon)\,e^{(2e^{1/2}d_2nt\epsilon)^2}-1\right)\end{equation}
\end{lem}
\proof
 For any $n\geq 1$ we have
\begin{align*}
\EE\left(\Vert \Ea_{s,t}(A^{\epsilon})\Vert^n\right)~\leq&~ \EE\left[
\exp{\left[n\int_s^t\Vert A^{\epsilon}_u\Vert\,du\right]}\right]\\
=&~1+\sum_{k\geq 1}\,\frac{n^k}{k!}\,\EE\left(\left[\int_s^t\Vert A^{\epsilon}_u\Vert\,du\right]^k\right)\\
\leq &~1+\sum_{k\geq 1}\,\frac{n^k}{k!}\,\left[\int_s^t\EE\left(\Vert A^{\epsilon}_u\Vert^k\right)^{1/k}\,du\right]^k \\
\leq&~ 1+\frac{1}{2}\,\sum_{k\geq 1}\,\frac{(2n(t-s))^k}{k!}\,(d_1^k+d_2^k\,\epsilon^k\,k^{k/2})
\end{align*}
This implies that 
\begin{eqnarray*}
\EE\left(\Vert \Ea_{s,s+t}(A^{\epsilon})\Vert^n\right)&\leq& \frac{1}{2}\,e^{2d_1nt}+\frac{1}{2}\sum_{k\geq 1}\,\frac{(2^{3/2}d_2nt\epsilon)^k}{k!}\,(k/2)^{k/2}
\end{eqnarray*}
On the other hand, by Stirling approximation we have
\begin{eqnarray*}
\frac{k^{k}}{(2k)!}&\leq&
\frac{k^{ k}}{\sqrt{4\pi k}\,(2k)^{2k}\,e^{-2k}} ~=~ \frac{e^{k}}{\sqrt{4\pi }\,2^{2k}} 
 \frac{1}{\sqrt{k}\,k^{k}\,e^{-k}} ~\leq~  \frac{e}{\sqrt{4\pi }}\,\frac{(e/4)^{k}}{k!} 
\end{eqnarray*}
and
\begin{eqnarray*}
\frac{(2k+1)^{k+1/2}}{(2k+1)!}&\leq&\frac{1}{\sqrt{2\pi(2k+1)}}\,\frac{1}{(2k+1)^{(2k+1)/2}\,e^{-(2k+1)}}\\
&\leq &\frac{e^2}{\sqrt{4\pi}}\,\left(\frac{e}{2}\right)^k\,\frac{1}{e\,\sqrt{k}\,k^{k}\,e^{-k}} ~\leq~ \frac{e^2}{\sqrt{4\pi}}\,\left(\frac{e}{2}\right)^k\,\frac{1}{k!}\\
\end{eqnarray*}
This implies that
$$
\sum_{k\geq 1}\,\frac{(2^{3/2}d_2n t\epsilon)^{2k}}{(2k)!}\,k^{k} ~\leq~ \frac{e}{\sqrt{4\pi }} \sum_{k\geq 1}\,\frac{((2e)^{1/2}d_2nt\epsilon)^{2k}}{k!} ~\leq~ \frac{e}{\sqrt{4\pi }}\,\left[e^{((2e)^{1/2}d_2n t\epsilon)^2}-1\right]
$$
and
\begin{align*}
&\sum_{k\geq 0}\,\frac{(2^{3/2}d_2nt\epsilon)^{2k+1}}{(2k+1)!}\,(2k+1)^{k+1/2}\\
&\qquad\qquad\leq~ (2^{3/2}d_2nt\epsilon)\left[1+\sum_{k\geq 1}\,\frac{(2^{3/2}d_2nt\epsilon)^{2k}}{(2k+1)!}\,(2k+1)^{k+1/2}\right]\\
&\qquad\qquad\leq~ (2^{3/2}d_2nt\epsilon)\left[1+\frac{e^2}{\sqrt{4\pi}}\sum_{k\geq 1}\,\frac{(2e^{1/2}d_2nt\epsilon)^{2k}}{k!}~\right]
~\leq~ 2^{1/2}\frac{e^2}{\sqrt{\pi}}\,d_2nt\epsilon\,e^{(2e^{1/2}d_2nt\epsilon)^2}
\end{align*}
This implies that
\begin{eqnarray*}
\EE\left(\Vert \Ea_{s,s+t}(A^{\epsilon})\Vert^n\right)&\leq& \frac{1}{2}\,e^{2d_1nt}+\frac{1}{4}\frac{e}{\sqrt{\pi}}
\left(e^{((2e)^{1/2}d_2n t\epsilon)^2}-1+2^{3/2}e\,d_2nt\epsilon\,e^{(2e^{1/2}d_2nt\epsilon)^2}\right)\\
&\leq&  \frac{1}{2}\,e^{2d_1nt}+\frac{e}{2}\,\sqrt{\frac{e}{\pi}}
\left(e^{((2e)^{1/2}d_2n t\epsilon)^2}-1+(2e)^{1/2}\,d_2nt\epsilon\,e^{(2e^{1/2}d_2nt\epsilon)^2}\right)
\end{eqnarray*}
We conclude that
\begin{eqnarray*}
\EE\left(\Vert  \Ea_{s,s+t}(A^{\epsilon})\Vert^n\right) &\leq&   \frac{1}{2}\,e^{2d_1nt}+\frac{e}{2}\,\sqrt{\frac{e}{\pi}}
\left((1+(2e)^{1/2}\,d_2nt\epsilon)\,e^{(2e^{1/2}d_2nt\epsilon)^2}-1\right)
\end{eqnarray*}
This ends the proof of (\ref{ref-nE-exp-Aa-1}). \qed

 Now we come to the proof of Theorem~\ref{theo-A-epsilon}. {\bf Proof of Theorem~\ref{theo-A-epsilon}:} Under $(H_1)$ we have the log-norm estimate
\begin{eqnarray*}
\frac{1}{t}\log{\Vert \Ea_{s,s+t}(\overline{A}^{\,\epsilon})\Vert} &\leq& \frac{1}{t} \int_0^t\,\mu(\overline{A}^{\,\epsilon}_{s+u})\,du\\
&\leq & \mu(A_{\infty})+\frac{1}{t} \int_0^t\,\Vert A_{\infty}-A_{s+u}\Vert\,du+
 	\frac{1}{t} \int_0^t\,\Vert A_{s+u}-\overline{A}^{\,\epsilon}_{s+u}\Vert\,du\\
&\leq& \mu(A_{\infty})+\frac{a}{b}\,\frac{e^{-as}}{t}+c_{2}\,\epsilon
\end{eqnarray*}
In the last assertion we have used the fact that 
$$
(H_1) \qquad\Longrightarrow\qquad \Vert A_{s+u}-\overline{A}^{\,\epsilon}_{s+u}\Vert \,\leq\, \EE(\Vert A_{s+u}-A^{\epsilon}_{s+u}\Vert^2)^{1/2} \,\leq\, c_{2}\,\epsilon
$$
Observe that
\begin{equation*}
e^{-as}\wedge \frac{1}{t}
\,\leq\, \frac{b}{4a}\vert\mu(A_{\infty})\vert\quad \mbox{\rm and}\quad \epsilon\leq \epsilon_2\wedge \frac{\vert\mu(P_{\infty})\vert}{4c_2}\quad\Longrightarrow\quad
\frac{1}{t}\log{\Vert \Ea_{s,s+t}(\overline{A}^{\,\epsilon})\Vert^2}\,\leq\, \mu(A_{\infty})
\end{equation*}
This completes the proof of the l.h.s. of (\ref{EA-epsilon-H1-}). Arguing as above we have the log-norm estimate
$$
\displaystyle\frac{1}{t}\log{\Vert \Ea_{s,s+t}(A^{\epsilon})\Vert}-\mu(A_{\infty}) \,\leq\, \frac{1}{t}\,\int_0^{t}\,\Vert A_{s+u}-A^{\epsilon}_{s+u}\Vert\,du+\frac{a}{b}\,\frac{e^{-as}}{t}
$$
 Observe that
 $$
s\vee t\geq h \quad\Longrightarrow\quad  e^{-as}\wedge \frac{1}{t}
 \,\leq\, \frac{b}{4a}\,\vert\mu(A_{\infty})\vert
$$
Taking the expectation we check the r.h.s. of (\ref{EA-epsilon-H1-}).

We consider the collection of events
$$
\Omega^{\epsilon}_{s,t}\,:=\,\left\{\Vert \Ea_{s,s+t}(A^{\epsilon})\Vert \,\leq\, \exp{\left[
 \frac{t}{2}\,\mu(A_{\infty})\right]}\right\}
 $$
 Assume $(H_1)$ holds. In this case, applying the Markov inequality, for any $\epsilon\leq \epsilon_n$ and any $s\vee t\geq h$
we have the uniform estimates
\begin{equation}
\begin{array}{l}
\displaystyle\EE\left(\Vert\frac{1}{t}~\int_0^{t}\,\Vert A_{s+u}-A^{\epsilon}_{s+u}\Vert\,du\Vert^n\right)^{1/n} \,\leq\,
\epsilon\,c_n\\
\\
\displaystyle \qquad\Longrightarrow\quad
\sup_{s\vee t\geq h}{\PP\left(\Omega-\Omega^{\epsilon}_{s,t}\right)}^{1/n} \,\leq\, \epsilon\,\overline{c}_n\quad\mbox{\rm with}\quad
\overline{c}_n:=4c_n/\vert\mu(A_{\infty})\vert
\end{array}
\label{def-overline-cn}
\end{equation}
Applying Fatou's Lemma for any $h\leq s \leq t_k\uparrow_{k\rightarrow\infty}\infty$,we have
$$
\begin{array}{l}
\displaystyle\PP\left(\forall l\geq 1~~\exists k\geq l~:~\Vert \Ea_{s,s+t_k}(A^{\epsilon})\Vert\,>\, \exp{\left[
 \frac{t_k}{2}\,\mu(A_{\infty})\right]}\right)\\
 \\
\qquad\qquad\displaystyle\leq~
\liminf_{k\rightarrow\infty}\,\PP\left(
\Vert \Ea_{s,s+t_k}(A^{\epsilon})\Vert\,>\, \exp{\left[
 \frac{t_k}{2}\,\mu(A_{\infty})\right]}\right)
 \leq \epsilon^n\,\overline{c}^n_n
 \end{array}
$$ 
This implies that for any $\epsilon\leq \epsilon_n$ we have 
$$
 \PP\left(\limsup_{k\rightarrow\infty}\,\frac{1}{t_k}\,\log{\Vert \Ea_{s,s+t_k}(A^{\epsilon})\Vert^2}\,<\, \mu(A_{\infty}) \right) \,\geq\,  1-\epsilon^n\,\overline{c}^n_n
$$
This ends the proof of (\ref{EA-epsilon-H1}).

Using (\ref{ref-nE-exp-Aa-1}) we have the rather crude estimate
\begin{eqnarray*}
d_2\,n\,t\epsilon^2\,\leq\, 1 \quad\Longrightarrow\quad
\EE\left(\Vert  \Ea_{s,s+t}(A^{\epsilon})\Vert^n\right)&\leq&   \frac{1}{2}\,e^{2d_1nt}+e\,\sqrt{\frac{e}{\pi}}
e^{(4e+(2e)^{1/2})d_2nt}\\
&\leq &2\,e^{dnt}\quad \mbox{\rm with}\quad d=(4e+(2e)^{1/2})\,(d_1\vee d_2)
\end{eqnarray*}
Using the Cauchy-Schwartz inequality for any $n\geq 2$ and any $\epsilon$ such that $d_2nt\epsilon^2\leq 1$ we have
\begin{eqnarray*}
\EE\left(\Vert \Ea_{s,s+t}(A^{\epsilon})\Vert^{n/2}\right)^{2/n}&=&\EE\left(\Vert \Ea_{s,s+t}(A^{\epsilon})\Vert^{n/2}~1_{\Omega^{\epsilon}_{s,t}}\right)^{2/n}+
\EE\left(\Vert\Ea_{s,s+t}(A^{\epsilon})\Vert^{n/2}~1_{\Omega-\Omega^{\epsilon}_{s,t}}\right)^{2/n}\\
&\leq & \exp{\left[
 \frac{t}{2}\,\mu(A_{\infty})\right]}+\epsilon\,\overline{c}_n\,2^{1/n}\,\exp{[td]}\\
 &\leq & (1+\overline{c}_n\,2^{1/n})\,\exp{\left[
 \frac{t}{2}\,\mu(A_{\infty})\right]}\leq \exp{\left[
 \frac{t}{4}\,\mu(A_{\infty})\right]}
 \end{eqnarray*}
 with the constant $\overline{c}_n$ defined in (\ref{def-overline-cn}) as soon as
 $$
\frac{4}{\vert \mu(A_{\infty})\vert}\,\log{\left(1+\overline{c}_n~2^{1/n}\right)} \,\leq\,  t \,\leq\, \frac{1}{d+\vert \mu(A_{\infty})\vert/2}\,\log{(1/\epsilon)}
 $$
 For instance, we can choose $\epsilon$ sufficiently small s.t. 
$$
c_{n}^- \,\leq\, t \,\leq\, c_{n}^+\log{(1/\epsilon)}
$$ 
with
$$
c_{n}^-:=\frac{4}{\vert \mu(A_{\infty})\vert}~\log{\left(1+\overline{c}_n~2^{1/n}\right)}\quad \mbox{\rm and}\quad
c_{n}^+:=\left(\frac{1}{d+\vert \mu(A_{\infty})\vert/2}\wedge \frac{1}{2d_2n} \right)
$$
In summary, for any $n\geq 1$ there exists some finite constants
$c_n^-\leq c^+_n$ such that for any $\epsilon\leq \exp{[-(c^-_n/c^+_n)}]$ and any $s\geq 0$ we have
$$
c_{n}^-\leq t\leq c_{n}^+\log{(1/\epsilon)} \quad\Longrightarrow\quad \frac{1}{nt}\log{\EE\left(\Vert \Ea_{s,s+t}(A^{\epsilon})\Vert^{n}\right)} \,\leq \,
 \frac{1}{4}~\mu(A_{\infty})
$$
This ends the proof of  (\ref{EA-epsilon-H2}). The proof of the theorem is now complete. \qed

\end{document}